\documentclass[10pt]{article}
\usepackage{amsmath,theorem,array}
\usepackage{tabularx,xspace,doublespace}
\usepackage{amssymb}
\usepackage{commands,shell}
\ifpdf
\usepackage[pdftex]{graphicx}
\else
\usepackage[dvips]{epsfig}
\fi
\begin{document}
\title{Output Feedback Control of Jet Engine Stall and Surge \\ 
Using
Pressure Measurements\thanks{This work 
was supported by NASA Glenn Research Center, Grant NAG3-2084.}}
\author{\bf
Manfredi Maggiore\hspace{1cm} 
Kevin Passino \\
\\
\it
Department of Electrical Engineering,
The Ohio State University \\
\it 
 2015 Neil Avenue, Columbus, OH
43210-1272 \\
}
\date{}
\maketitle
\begin{abstract}
The problem of controlling surge  and stall in jet engine compressors
is of fundamental importance in preventing damage and lengthening the
life of these components. In this paper, we use the Moore-Greitzer 
mathematical model to develop an output feedback controller for  these two instabilities
(only one of the three states is measurable). 
This problem is particularly challenging since the system is not
completely observable and, hence, none of the output feedback control techniques found in
the literature can be applied to recover the performance of a full
state feedback controller. However, we show how to successfully
solve it by using a novel output feedback approach for the
stabilization of general stabilizable and incompletely observable
systems.  
\end{abstract}
\vspace*{-1cm}
\section{Introduction and Problem Description}
In this paper  we consider the problem of controlling two 
instabilities  which occur in jet engine compressors, namely 
rotating stall and surge. Rotating stall develops when there is a
region of stagnant flow rotating around the circumference of the
compressor causing undesired vibrations in the blades
and reduced pressure rise of the compressor. 
Surge is an axisymmetric oscillation of the flow through the
compressor that can cause undesired vibrations in other components of
the compression system and damage to the engine.
In \cite{MooGre86}, Moore and Greitzer developed a three-state finite
dimensional Galerkin approximation of a nonlinear PDE model
describing the compression system. Since its development, several
researchers have used the Moore-Greitzer three state model (MG3) to
design stabilizing controllers for stall and surge. The available
control approaches may be divided into three main
categories: 1) Linearization and  linear perturbation models (e.g.,
\cite{WilHeeJagStoo99,PadValEpsGreGue94,EpsWilGre89} among others); 2) Bifurcation
analysis (e.g.,
\cite{LiaAbe92,LiaAbe96,EveGysNetSha98,McC90,AbeHouHos93}); and  3)
Lyapunov based methods (e.g., \cite{KrsFonKokPad98,BadChoNet96,WanKrsLar97}).
Most  existing results focus on the development of state
feedback controllers, thus complicating  their practical
implementation as in \cite{EveGysNetSha98}, where the authors
use sensor arrays (2D sensing) to implement a state feedback control
law depending on the squared  amplitude of the first harmonic of
asymmetric flow and the derivative of the  air flow through the
compressor. In \cite{KrsFonKokPad98}, a partial
state feedback controller simplifies practical implementation by only
requiring measurements of the mass flow and plenum pressure rise
(hence 2D sensing is not needed).
 On the other hand, the limitation
of this partial
state feedback controller lies in the fact that it 
cannot globally stabilize  a unique equilibrium point. 

To the best of our knowledge, no attempt has been made to design 
a stabilizing output feedback controller (using only plenum
pressure rise feedback) based on a full-state
feedback control law.  This is probably due to the fact that MG3
becomes unobservable when there is no mass flow through the
compressor, \ie the system is not uniformly completely
observable (UCO), and none of the techniques found in the output
feedback control literature (e.g., \cite{EsfKha92,Tor92,TeePra94, KhaEsf93,MahKha96,MahKha97,AtaKha99}) can be employed for the solution of this
problem. 
In this paper we introduce
a new globally stabilizing full state feedback control law for MG3, and
we employ the theory developed in \cite{MagPas00-1,MagPas00-3} for the
output feedback control of incompletely observable nonlinear systems to
regulate  stall and surge by using only pressure measurements.
The MG3 model  is described by
(see \cite{KrsKanKok95,KrsFonKokPad98} for an analogous exposition)
\begin{equation}
\begin{split}
&\dot \Phi = -\Psi + \Psi_C ( \Phi ) - 3 \Phi R \\
&\dot \Psi = \frac{1}{\beta^2} ( \Phi - \Phi_T ) \\
&\dot R = \sigma R (1-\Phi^2-R), \, \, R(0)\geq 0
\label{eq:sys:1} 
\end{split}
\end{equation}
where $\Phi$ represents the mass flow, $\Psi$ is the plenum pressure rise, $R
\geq 0$ is the normalized stall cell squared amplitude, $\Phi_T$ is the mass
flow through the throttle,  $\sigma = 7$,  and $\beta = 1/\sqrt{2}$.
The functions $\Psi_c(\Phi)$ and $\Phi_T(\Psi)$ are the compressor and
throttle characteristics, respectively, and are defined as
$
\Psi_C(\Phi)= \Psi_{C_0}+1+3/2 \Phi - 1/2 \Phi^3$,
$\Psi =  \frac{1}{\gamma} (1 + \Phi_T(\Psi))^2$,
where $\Psi_{C_0}$ is a constant and $\gamma$ is the throttle opening,
the control input. Given the static relationship existing between
$\Phi_T$ and $\gamma$, without loss of generality, in what follows we
will design a controller assuming that $\Phi_T$ is our control input.
Our control objective is to stabilize system (\ref{eq:sys:1})
around the critical equilibrium $R^e =  0 , \Phi^e  = 1 , \Psi^e =
\Psi_C(\Phi^e) = \Psi_{C_0} + 2$, which achieves the peak operation on
the compressor characteristic.
We shift the origin to the
desired equilibrium with the change of variables $\phi = \Phi - 1 ,
\psi = \Psi - \Psi_{C_0} -2$. System (\ref{eq:sys:1}) then becomes
\begin{equation}
\begin{split}
&\dot R = -\sigma R^2 - \sigma R ( 2 \phi + \phi^2 ) \\
&\dot \phi = -\psi - 3/2 \phi^2 - 1/2 \phi^3 - 3 R \phi - 3 R \\
&\dot \psi = - \frac{1}{\beta^2}(\Phi_T - 1 - \phi)
\label{eq:sys}
\end{split}
\end{equation}
The pressure rise (and hence $\psi$) is the only measurable state variable.
It is readily seen that this system is input output \FLE
with relative degree one (the first derivative of  
$\psi$ contains the input $\gamma$), and its zero-dynamics are
nonminimum phase.
\section{State Feedback Control Design}
For convenience,
in the remainder of the paper we will 
redefine the control input to be $u= \Phi_T - 1$. 
Next, notice that
Assumption A2 in \cite{MagPas00-1,MagPas00-3} is satisfied since, for example, a stabilizing control law for 
 (\ref{eq:sys}) is given in \cite{KrsKanKok95} by means of
backstepping design. However, the control law proposed in
\cite{KrsKanKok95} turns out to be quite complex.
In \cite{KrsFonKokPad98},  it is shown that a linear partial state
feedback control law of the type $u= d_1 \psi - d_2 \phi$ achieves
either a unique asymptotically stable equilibrium point with domain of
attraction $\{ (R,\phi,\psi) \in \Re^3| R \geq 0 \}$ or two equilibria
on the axisymmetric and stall characteristic, with domains of
attraction $\{ (R,\phi,\psi) \in \Re^3| R = 0 \}$ and $\{
(R,\phi,\psi) \in \Re^3| R > 0 \}$, respectively (see Theorem 3.1 in
\cite{KrsFonKokPad98}). 
Here, this problem is overcome by viewing system
(\ref{eq:sys}) as an interconnection of two subsystems, namely the
$R$-subsystem and the $(\phi, \psi)$-subsystem, and then building  a full
 state feedback controller which makes the origin of \refeq{eq:sys} an
asymptotically stable equilibrium point with domain of attraction $\{
(R,\phi,\psi) \in \Re^3| R \geq 0 \}$,  as seen in
the next theorem. 
\begin{thm}
For system (\ref{eq:sys}), with the choice of the control law
\begin{equation}
\bar u = (1 - \beta^2 k_1 k_2)  \phi + \beta^2 k_2 \psi + 3 \beta^2
 k_1 R \phi
\label{eq:control} 
\end{equation}
where $k_1$ and $ k_2$ are positive scalars satisfying the inequalities,
\begin{align}
&k_1 > \frac{17}{8} + \frac{(2 C \sigma + 3)^2}{2} \label{eq:k1_1} \\
&\left(C \sigma - \frac{105}{64} \right) k_1^2 + \frac 3 4 \left(
-\frac 1 2 C \sigma + \frac {21} 4 \right) k_1 - (C \sigma + 3)^2> 0
\label{eq:k1_2} \\
&k_2 > k_1 + \frac{9}{4} k_1^2 + \frac{9 k_1}{4 k_1 - 9/2} +
\frac{(k_1^2 -1)^2}{4} \label{eq:k2} \\
& C > \frac{3}{2 \sigma}\label{eq:C}
\end{align}
the origin is an asymptotically stable equilibrium point with domain of
attraction $\mathcal{A} = \{ (R,\phi,\psi) \in \Re^3| R \geq 0 \}$.
\end{thm}
\pf
For the sake of simplicity, redefine the control input to be
$u' = -\frac{1}{\beta^2}(u - \phi)$,
so that the last equation in (\ref{eq:sys}) becomes
$\dot \psi = u'$. Next, notice that system (\ref{eq:sys}) can be
viewed as the interconnection of two subsystems:

\vspace*{-1.3cm}
\begin{equation}
\begin{aligned}[cl]
[S_1] \,\,
 \dot R = -\sigma R^2 , \s &
[S_2]   \left\{ \begin{aligned}[c]
\dot \phi &= -\psi - \frac{3}{2} \phi^2 - \frac{1}{2} \phi^3 \\
\dot \psi &= -u' \nonumber
\end{aligned}\right.
\end{aligned}
\end{equation}
A Lyapunov function for $[S_1]$, defined on the domain $\{ R \in \Re \,|\, R
\geq 0 \}$, is $V_1 = R$, and its time derivative is readily
found to be
$
\dot V_1 = - \sigma R^2 
$ 
thus showing that the origin of $[S_1]$ is an \asy stable equilibrium
point of $[S_1]$, and its domain of attraction is $\{ R \in \Re \,|\, R
\geq 0 \}$.
As for subsystem $[S_2]$ the analysis found in Section 2.4.3 in \cite{KrsKanKok95}
suggests using $V_2 = \frac{1}{2} \phi^2 + \frac{k_1}{8} \phi^4 + \frac{1}{2}(\phi -
k_1 \psi)^2$,
where $k_1$ is a positive design constant. Furthermore, in
\cite{KrsKanKok95}, a stabilizing control law for $[S_2]$ is found to be
$u' = -c_1 \phi + c_2 \psi$, where $c_1$ and $c_2$ are two appropriate
positive constants. In the following we will show that, in order to
stabilize the interconnection of systems $[S_1]$ and $[S_2]$, one needs to
add to $u' = -c_1 \phi + c_2 \psi$ a term which is proportional to the product $R
\phi$.
Based on these considerations, consider the following candidate
Lyapunov function for system (\ref{eq:sys}),
\begin{equation}
V = C V_1 + V_2 =  C R + \frac{1}{2} \phi^2  + \frac{k_1}{8} \phi^4 + \frac{1}{2} 
\left( \psi - k_1 \phi \right)^2 
\label{eq:V}
\end{equation}
where $C >0$ is a scalar.
After noticing that $V$ is positive definite on the domain
$\mathcal{A}$, and letting $\tilde \psi = \psi - k_1 \phi$, we
calculate the time derivative of V as follows,

\begin{equation}\begin{split}
\dot V = &- C \sigma R^2 - C \sigma R (2 \phi + \phi^2) + \left(\phi +
\frac{k_1}{2} \phi^3 \right) \left(-\psi -\frac{3}{2} \phi^2 -
\frac{1}{2} \phi^3 - 3 R \phi - 3 R \right) + \\
&+ \tilde \psi \left( u' +
k_1 \psi + \frac{3}{2} k_1 \phi^2 + \frac{1}{2} k_1 \phi^3 + 3 k_1 R
\phi + 3 k_1 R \right) 
\end{split}
\label{eq:Vdot}
\end{equation}
Here, as in  \cite{KrsKanKok95}, we use the identity 
$
-\frac{3}{2} \phi^2 - \frac{1}{2} \phi^3 = -\frac{1}{2} \left(\phi +
 \frac{3}{2}\right)^2 \phi + \frac{9}{8} \phi
$
to  eliminate the potentially  destabilizing term $-\left( \phi +
k_1/2 \phi^3 \right) 3/2 \phi^2$.
Next,  substituting  (\ref{eq:control}) into (\ref{eq:Vdot}) (after taking
in account the definition of $u'$), letting
$\bar k_1 = k_1 - 9/8$, and using 
the definition of $\tilde \psi$, we get
\begin{equation}\begin{split}
\dot V = & - C \sigma R^2 - C \sigma R (2 \phi + \phi^2) + \left(\phi +
\frac{k_1}{2} \phi^3 \right) \left(-\tilde \psi - \bar k_1 \phi -\frac{1}{2} \left(\phi +
 \frac{3}{2}\right)^2 \phi 
 - 3 R \phi - 3 R \right) + \\
&+ \tilde \psi \left( -(k_2 - k_1) \tilde \psi + k_1^2  \phi 
 + \frac{3}{2} k_1 \phi^2 + \frac{1}{2} k_1 \phi^3  + 3 k_1 R \right)
\end{split}\end{equation}
Now notice that   the expression  $-\left(\phi +
\frac{k_1}{2} \phi^3 \right) \frac{1}{2} \left(\phi +
 \frac{3}{2}\right)^2$ can be discarded since it  is negative
definite, and that the term $\frac{k_1}{2}  \phi^3 \tilde \psi$ cancels
out. After collecting the remaining  terms, we  get 
\begin{equation}\begin{split}
\dot V \leq & - C \sigma R^2 - (2 C \sigma + 3) R \phi - (C \sigma +
3) R \phi^2 - \bar k_1 \phi^2  - \left(\frac{k_1 \bar k_1}{2} +
\frac{3 k_1}{2} R \right) \phi^4 - \frac{3 k_1}{2} R \phi^3 + \\
& +\tilde
\psi  \left( -(k_2 - k_1) \tilde \psi + (k_1^2 -1) \phi 
 + \frac{3}{2} k_1 \phi^2   + 3 k_1 R \right)
\label{eq:Vdot2}
\end{split}\end{equation}
By using Young's inequality five times we have
\begin{gather*}
-(2 C \sigma + 3) R \phi \leq \frac{1}{2} R^2 + \frac{(2 C \sigma +
3)^2}{2} \phi^2 ,
\s -\frac{3 k_1}{2} R \phi^3 \leq \frac{3 k_1}{2} \left( \frac{R
\phi^2}{4} + R \phi^4 \right) , \\
(k_1^2 -1) \phi \tilde \psi \leq \phi^2 + \frac{(k_1^2 - 1)^2}{4}
\tilde \psi^2 ,
\s 3  k_1  R \tilde \psi  \leq R^2 + \frac{9}{4} k_1^2 \tilde
\psi^2 ,
\s \frac{3}{2} k_1 \phi^2  \tilde \psi \leq \frac{k_1 \bar k_1}{4}
\phi^4 + \frac{9 k_1}{4 \bar k_1} \tilde \psi^2
\end{gather*}
Applying the inequalities above to (\ref{eq:Vdot2}) we get
\begin{align}
\dot V \leq & - \left( C \sigma - \frac{3}{2}\right) R^2 - \left(\bar k_1
- \frac{(2 C \sigma +3)^2}{2} -1 \right) \phi^2 - \left(k_2 - k_1 -
\frac{9}{4} k_1^2 - \frac{9 k_1}{4 \bar k_1} - \frac{(k_1^2 -1)^2}{4}
\right) \tilde \psi^2 + \nonumber \\
& - \left(C \sigma +3 - \frac{3}{8} k_1 \right) R \phi^2 - 
\frac{k_1 \bar k_1}{4} \phi^4 , \nonumber \\
\leq & -\begin{bmatrix}
R \\ \phi^2
\end{bmatrix}^\top 
\begin{bmatrix}
C \sigma - \frac{3}{2} & \frac{1}{2}\left(C \sigma +3 - \frac{3}{8}
k_1 \right)\\
\frac{1}{2}\left(C \sigma +3 - \frac{3}{8} k_1 \right) & \frac{1}{4}k_1 \bar k_1
\end{bmatrix}
\begin{bmatrix}
R \\ \phi^2
\end{bmatrix} - \left(\bar k_1
- \frac{(2 C \sigma +3)^2}{2} -1 \right) \phi^2  +  \nonumber\\
& - \left(k_2 - k_1 -
\frac{9}{4} k_1^2 - \frac{9 k_1}{4 \bar k_1} - \frac{(k_1^2 -1)^2}{4}
\right) \tilde \psi^2 
\end{align}
Hence,  $\dot V$ is negative definite on the domain $\mathcal{A}$,
provided that the quadratic form above is positive definite and that
the coefficients multiplying $\phi^2$ and $\tilde \psi^2$ be positive.
By imposing the positive definiteness of the quadratic form we obtain
$C \sigma - \frac{3}{2} > 0$,
$\left(C \sigma - \frac{3}{2}\right) \frac{1}{4}k_1 \bar k_1 -
\frac{1}{4} \left(C \sigma +3 - \frac{3}{8} k_1 \right)^2 > 0$,
while by imposing the positivity of the coefficients of the remaining
two terms we get
$
\bar k_1 > \frac{(2 C \sigma +3)^2}{2} + 1$,
 $k_2 > k_1 + \frac{9}{4} k_1^2 + \frac{9 k_1}{4 \bar k_1} + \frac{(k_1^2 -1)^2}{4}$.
By using the definition of $\bar k_1$, inequalities (\ref{eq:k1_1}),
\refeq{eq:k1_2}, (\ref{eq:k2}), and (\ref{eq:C}) follow.
 In conclusion, if $k_1$, $k_2$, and $C$
are chosen so that (\ref{eq:k1_1})-(\ref{eq:C}) hold, we have that
$\dot V$ is negative definite on $\mathcal{A}$ which contains the origin. 
This  leads to the
conclusion  that $\{ R =0, \phi=0, \tilde \psi = 0 \}$ is an \asy stable
equilibrium point, which in turn implies that $\{ R=0, \phi=0, \psi=0
\}$ is an \asy stable equilibrium point.
Our next objective is to show that $\mathcal{A}$ is a region of
attraction for the origin. This, however, is not immediately evident
from our result, since the set $\{ [ R, \phi, \psi ]^\top \in \Re^3
\,|\, V \leq K, K>0 \}$ is unbounded and, due to the presence of the term
$C R$ in $V$,  it is not completely contained
in $\mathcal{A}$. In other words, it may happen that, while the Lyapunov
function is decreasing, $R$ becomes negative, and thus the state
trajectory exits the set  $\mathcal{A}$, where $\dot V$ is guaranteed
to be negative definite. Therefore, in order to complete our analysis,
we need to show that $\mathcal{A}$ is invariant, which, together with
$\dot V < 0$, implies that
the set $\{ [ R, \phi, \psi ]^\top \in \Re^3 \,|\, V \leq K, K>0 \}
\cap \mathcal{A}$ is a region of attraction of the origin for any $K
>0$. This is readily seen by noticing that, on the boundary of
$\mathcal{A}$, $R =0$. From (\ref{eq:sys}), $R=0$ implies $\dot R =
0$, thus proving that no trajectory of the system can cross the
boundary of $\mathcal{A}$, and  therefore  $\mathcal{A}$ is invariant.
In conclusion, given any initial condition $[ R(0), \phi(0), \psi(0)
]^\top$ in $\mathcal{A}$, there exists a constant $K > 0$ such that the
initial condition is contained in the set  $\{ [ R, \phi, \psi ]^\top 
\in \Re^3 \,|\, V \leq K, K>0 \} \cap \mathcal{A}$,  thus proving that
the origin  of system (\ref{eq:sys}) is an \asy stable equilibrium point
with domain  of attraction $\mathcal{A}$.
\BBOX
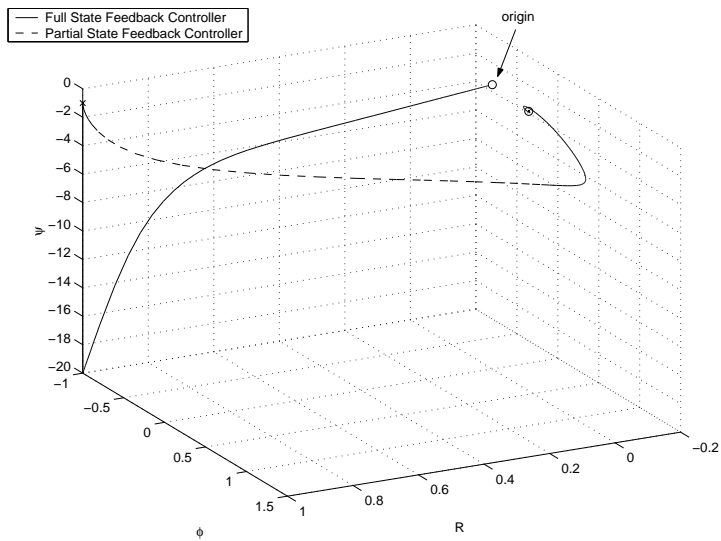
\scalefig{sf}{.6}{Comparison between the
partial state feedback controller developed in \cite{KrsFonKokPad98}
and the full state feedback  controller \refeq{eq:control}.} 
\begin{rem}
By using inequalities \refeq{eq:k1_1}-\refeq{eq:C},
it is easy to show that the only equilibrium point of the \cls on the
set $\mathcal{A}$ is the origin, as predicted by Theorem 1. \refig{sf} shows the evolution of
the closed-loop trajectories under  the
partial state feedback controller developed in \cite{KrsFonKokPad98}
and the controller \refeq{eq:control} for a particular choice of the
coefficients $d_1, d_2, k_1, k_2$. The partial state
feedback controller stabilizes an equilibrium point different from the
origin $(R, \phi, \psi) = (0,0,0)$. 
\end{rem}
\begin{rem}
Inequalities \refeq{eq:k1_1}-\refeq{eq:C} represent conservative
bounds on $k_1$ and $k_2$. In practical implementation, these
parameters may be chosen significantly smaller after some tuning.
\end{rem}
In order to complete the state feedback design, we have to add an
appropriate number of integrators at the input side of the system (see
\cite{MagPas00-1,MagPas00-3}).
Following the procedure outlined in  \cite{MagPas00-1,MagPas00-3}, we form the
observability mapping $\HH$
\begin{equation}
 y_e = 
\begin{bmatrix}
y \\
\dot y \\
\ddot y  
\end{bmatrix} = \HH\left([R, \phi, \psi]^\top, u, \dot u\right) =\begin{bmatrix}
\psi \\
- 1/ \beta^2 (u  - \phi) \\
1/\beta^2 \left(- \dot u -\psi - 3/2 \phi^2 - 1/2 \phi^3 - 3 R
\phi - 3 R \right)
\end{bmatrix}
\label{eq:observability:example}
\end{equation}
Notice that the observability assumption A1 in
\cite{MagPas00-1,MagPas00-3} is satisfied, for all $\phi \neq -1$, 
with $n_u = 2$ in that given $y_e , u$, and
$\dot u$, one can uniquely find $R, \phi, \psi$. The operating
point $\phi = -1$ corresponds to $\Phi = 0$, \ie no mass flow
through the compressor which is a condition we would like to avoid
during normal engine operation.
Since $n_u=2$, we extend the system with two integrators
$\dot z_1 = z_2$, $\dot z_2 = v$, $u = z_1$. 
To simplify the notation in the following, define
$x = [ R,  \phi,  \psi ]^\top$, and rewrite (\ref{eq:sys})
as $\dot x = f(x) + g(x) z_1$.  
Next, we find a stabilizing control law for the extended system by using 
the integrator backstepping lemma:
$ v = \dot{\alpha} - \tilde z_1 - k_4 \tilde z_2 \triangleq
\varphi(x,z)$,  where $\tilde z_1 = z_1 - \bar u$, 
$\alpha = -k_3 \tilde z_1 - \DER{V}{x} g(x) + \DER{\bar u}{x} [ f(x) + g(x) \,
z_1 ]$, $\tilde z_2 = z_2 - \alpha$, 
and  $k_3 , k_4 $ are arbitrary positive constants. 
This completes the design of a stabilizing state feedback for the
extended system. The Lyapunov function of the closed-loop extended
system is $
\bar V = V + \frac{1}{2} \tilde z_1^2 + \frac{1}{2} \tilde z_2^2
$.
Notice that, following the same reasoning as in  the proof of Theorem 1 the set 
$\{ [R, \phi, \psi, z_1, z_2]^\top \in \Re^5 \,|\, R \geq 0 \}$ is
invariant; hence by applying the backstepping lemma we guarantee that
the origin of the extended system is asymptotically stable with domain
of attraction $\DD = \mathcal{A} \times \Re^2$.
\vspace*{-.5cm}
\section{Output Feedback Design}
The validity of the observability assumption A1 in
\cite{MagPas00-1,MagPas00-3} allows us to design a stable observer.
As already pointed out, Assumption A1 in \cite{MagPas00-1,MagPas00-3} 
is satisfied on the domain $\X \times \U  = \left\{ [R , \phi, \psi] 
\in \Re^3 \,|\, \phi > -1 \right\} \times \Re^2$.
We first design the observer developed in \cite{MagPas00-1,MagPas00-3},
\begin{equation}
\begin{split}
&\dot {\hat R} = -\sigma {\hat R}^2 - \sigma R (2 \hat \phi + {\hat
\phi}^2 ) - \frac{(l_1 / \rho) + \beta^2 (3 \hat \phi + 3 \hat R + (3/2) {\hat
\phi}^2) (l_2/\rho^2)+ \beta^2 (l_3 / \rho^3)}{3 ( 1 + \hat \phi)} (\psi - \hat \psi) \\
&\dot{\hat  \phi} = - \hat \psi - 3/2 \,{\hat \phi}^2 - 1/2 \, {\hat \phi}^3 - 3
{\hat R} {\hat \phi} - 3 \hat R + \beta^2 (l_2 / \rho^2)   (\psi - \hat \psi) \\
&\dot{\hat \psi} = -\frac{z_1  - \hat \phi}{\beta^2} + (l_1 / \rho)  (\psi - \hat \psi) 
\end{split}
\label{eq:observer}
\end{equation}
where $\rho$ is a positive design parameter and the vector $L = [l_1,
l_2, l_3 ]^\top  \in \Re^3$ is chosen to be Hurwitz. Next, we
calculate the solution $P$ of the Lyapunov equation $P (A_c - L C_c) +
(A_c - L C_c)^\top P = -I$, where $(A_c, C_c)$ is a canonical
observable pair. In order to confine the observer estimates
to within the observable space, we implement the following projection,
\begin{equation*}
\begin{split}
&\xPdot = \left[ \DER{\HH}{\hat x} \right]^{-1} \left\{ 
\PP \left(\hat \xi, \dot{\hat \xi}, z, \dot z \right) - \DER{\HH}{z}
\dot z \right\} \\
&\PP(\hat \xi,\dot{\hat \xi},z, \dot z) = 
\begin{cases}
\D \dot{\hat \xi} - \Gamma \frac{N(\hat
\xi) \left(
N(\hat \xi,z)^\top \dot{\hat \xi} + N_z(\hat \xi, z)^\top \dot z
\right)} {N(\hat \xi,z)^\top \Gamma N(\hat \xi,z)} 
& \text{if }  N(\hat \xi,z)^\top \dot{\hat \xi} + N_z(\hat \xi,
z)^\top \dot z \geq 0 \text{ and }
\hat \xi \in \partial C_\xi(z)\\
\dot{\hat \xi} & \text{ otherwise}
\end{cases}
\end{split}
\end{equation*}

\vspace*{-1cm}
\noindent
where $\Gamma = ( S \E' )^{-1} ( S \E' )^{-1} $, 
$S = S^\top$ denotes the matrix square root of $P$, $\hat \xi =
\HH(\hat x,z)$, $\dot{\hat \xi} = \left\{ \DER{\HH}{\hat x} \dot{\hat
x} + \DER{\HH}{z} \dot z \right\}$, and $C_\xi(z)$ is the cube
\[
C_\xi(z)= \hspace*{-.1cm}\left\{ \xi \in \Re^3 \,|\, \xi_1 \in [a_1,b_1], \xi_2 \in \left[
-\frac{1}{\beta^2} (z_1 +a_2), -\frac{1}{\beta^2} (z_1 - b_2)\right],
\xi_3 \in \left[ \frac{1}{\beta^2} (-z_2 
- a_3), \frac{1}{\beta^2} (-z_2 + b_3)\right] \right\}
\]
which, when $a_2 < 1$, is contained in $\HH(\X,z)$, for all
$z$ (the scalars $a_i$, $b_i$, $i=1, 2, 3$ have to be chosen to
satisfy Assumption A3 in \cite{MagPas00-1,MagPas00-3}). 
Finally, $N(\hat \xi,z)$ and  $N_z(\hat \xi,z)$ are the normal vectors to the
boundary of $C_\xi (z)$ \wrt $\xi$ and  $z$, respectively, and are
given by
\begin{align*}
&N(\hat \xi, z) = \begin{cases}
[1,0,0]^\top \text{ if } \hat \xi_1 = b_1 & [-1, 0,0]^\top  \text{
if } \hat \xi_1 = a_1 \\
{}[0,1,0]^\top  \text{ if } \hat \xi_2 =  -\frac{1}{\beta^2} (z_1 -
b_2) & [0,-1,0]^\top \text{ if }  \hat \xi_2 =  -\frac{1}{\beta^2} (z_1
+ a_2) \\
{}[0,0,1]^\top \text{ if } \hat \xi_3 = \frac{1}{\beta^2} (-z_2 + b_3) &
[0,0,-1]^\top \text{ if } \hat \xi_3 = \frac{1}{\beta^2} (-z_2 - a_3)
\end{cases} 
\end{align*}
\begin{align*}
&N_z(\hat \xi, z) = \begin{cases}
[0,0]^\top \text{ if } \hat \xi_1 = b_1 \text{ or }  \hat \xi_1 = a_1
\\
{}\left[\frac{1}{\beta^2}, 0 \right]^\top \text{ if } \hat \xi_2 =  -\frac{1}{\beta^2} (z_1 -
b_2) & \left[-\frac{1}{\beta^2}, 0 \right]^\top \text{ if }  \hat
\xi_2 =  -\frac{1}{\beta^2} (z_1 +  a_2) \\
{}\left[0,\frac{1}{\beta^2}\right]^\top \text{ if } \hat \xi_3 =
\frac{1}{\beta^2}( -z_2 + b_3 )& \left[0,-\frac{1}{\beta^2}\right]^\top
\text{ if } \hat \xi_3 = \frac{1}{\beta^2} (-z_2 - a_3 )
\end{cases}
\end{align*}
Thus, the output feedback controller design is completed by 
letting $\hat v = \varphi (\xP,z)$, and Theorem 2 in
\cite{MagPas00-1,MagPas00-3} guarantees that the origin of the \cls,
controlled by $\hat v$, is \asy stable with domain of attraction
$\DD' \times \Omega_{c_2}^x$, where  $ \Omega_{c_2}^x \eqdef \{ [R,
\phi, \psi]^\top \,|\, \bar V \leq c_2, \text{ and } R \geq 0 \}$,
$c_2 >0$ is the largest scalar such that $\Omega_{c_2}^x \subset \{[R,
\phi, \psi]^\top \in \Re^3 \,|\, \phi > -1 \}$, and $\DD' \subset
\Omega_{c_2}^x$ can be made arbitrarily close to $ \Omega_{c_2}^x $ by
choosing $\rho$ in \refeq{eq:observer} small enough (see Theorem 2 in
\cite{MagPas00-1,MagPas00-3}). 
\section{Simulation Results}
Here we present the simulation results when the output
feedback controller developed in the previous section is applied
to system (\ref{eq:sys}). We choose 
 $k_1 = 25$ and $ k_2 = 1.1 \cdot10^5$ to fulfill inequalities 
\refeq{eq:k1_1}-\refeq{eq:C} in Theorem 1.
In order to choose the size of the compact set $C_\xi(z)$ so that
Assumption A3 in \cite{MagPas00-1,MagPas00-3} is satisfied, we may use the Lyapunov function $\bar V$
to calculate $\Omega_{c_2}^x$, choose $c_2$ small enough to guarantee
that $\Omega_{c_2}^x \subset \X$, and use $\HH$ to calculate bounds on
$\xi$ when $x \in \Omega_{c_2}^x$. However, a more practical way to
address the design of $C_\xi(z)$ consists of running a number of
simulations for the \cls under state feedback corresponding to several
initial conditions $[R(0), \phi(0), \psi(0)]^\top$, and calculating
upper and lower bounds for $\psi$, $\phi$, and $-\psi-3/2 \phi^2 - 1/2
\phi^3 - 3 R \phi - 3 R$: these will provide the values of $a_i, b_i$,
$i=1, 2, 3$, respectively.
By doing that, we found that whenever $[R(0),
\phi(0), \psi(0) ]^\top \in \Omega_0 \eqdef \{ [R,
\phi, \psi ]^\top \in \Re^3 \,|\, R \in [0, 0.1], \phi \in
[-0.1, 0.1], \psi \in [-0.5, 0.5] \}$,  we have that $a_1 = -2$,
$b_1=1$, $a_2 = -0.5$, $b_2=1$, $a_3=-0.5$, $b_3=0.3$ satisfy
Assumption A3 in \cite{MagPas00-1,MagPas00-3}.  We must point out
that our choice of $\Omega_0$ is rather conservative and is made
primarily for the sake of illustration. The actual domain of
attraction $\DD'$ under output feedback control is larger that
$\Omega_0$.
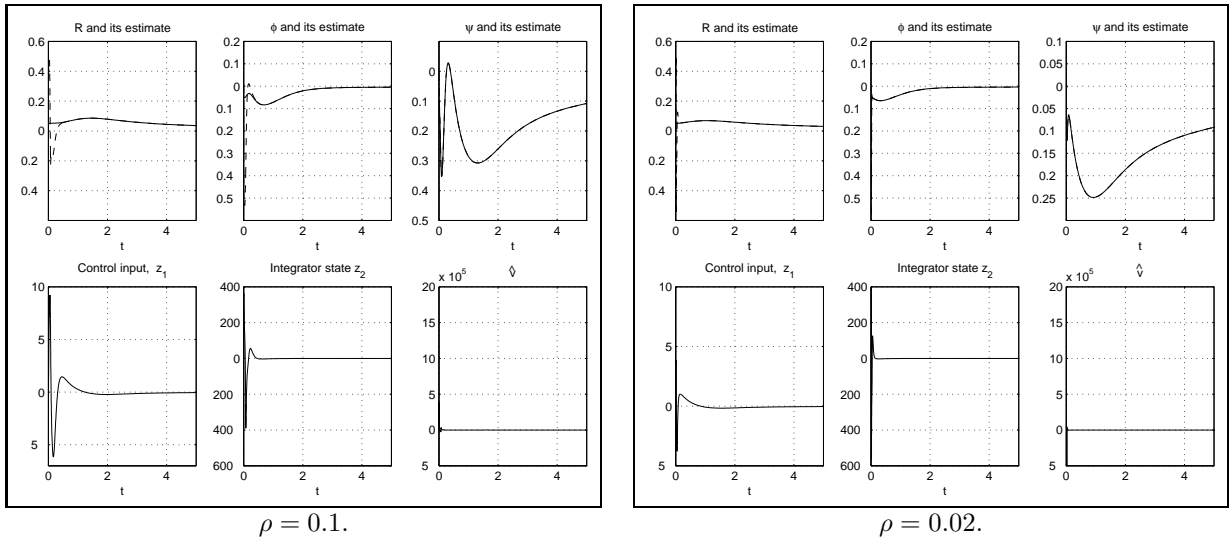
\twofigcap{simu_10}{$\rho =0.1$.}{simu_50}{$\rho=0.02$.}{Output feedback trajectories.}{.48}
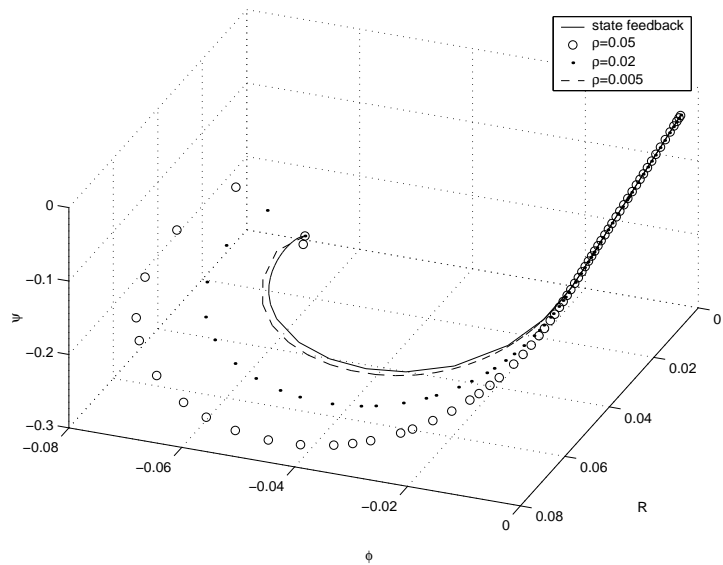
\scalefig{3D}{.6}{State feedback trajectories and output feedback trajectories
for $\rho=0.05$, $\rho=0.02$, and $\rho=0.005$.}
In Figure \ref{simu_10}  system and controller states, together with the
control input, are plotted for two decreasing values of $\rho$
confirming the theoretical predictions about the arbitrary fast rate of convergence
of the observer found in Theorem 1 in \cite{MagPas00-1}. Furthermore, the figures
 also show the operation of the
projection which prevents the observer from peaking and guarantees that $\hat
\phi> -0.5$. 
Finally, note that the output feedback
trajectories approach the state feedback ones, as 
showed in Figure \ref{3D}.
\nocite{MagPas00-2}
\bibliographystyle{abbrv}
\bibliography{biblio}
\end{document}